\documentclass[12pt]{article}
\usepackage{amssymb}
\usepackage{amsmath}
\newcommand{\la}{\lambda}
\newcommand{\lap}{\mbox{$\bigtriangleup$}}

\newcommand{\ra}{{\mbox{$\rightarrow$}}}
\newcommand{\be}{\begin{equation}}
\newcommand{\ee}{\end{equation}}

\newtheorem{mthm}{Theorem}

\newtheorem{mlem}{Lemma}

\newtheorem{thm}{Theorem}[section]

\newtheorem{lem}{Lemma}[section]

\begin{document}

\title{A Liouville Theorem for the
Higher Order Fractional Laplacian}

\author{ Ran Zhuo, \;Yan Li}
\date{\today}
\maketitle

\begin{abstract}
We deal with the higher-order fractional Laplacians by two methods: the integral method and the system method. The former depends on the  integral equation equivalent to the differential equation.
 The latter works directly on the differential equations.
We first
derive nonexistence of
 positive solutions, often known as the Liouville type theorem, for the integral and differential equations. Then through an delicate iteration, we show symmetry for positive solutions.

\end{abstract}


\section{Introduction}

Let $t=m+\frac{\alpha}{2}$ and $0<\alpha<2$. The higher-order
 fractional Laplacian $(-\lap)^{t}$ with $m \in  N$, and the usual
  fractional Laplacian $(-\lap)^{\alpha/2}$, usually take the form of
\begin{eqnarray*}
(-\lap)^{t}u(x)&=&C_{n,t}P.V.
\int_{R^n}\frac{\mathop\Sigma\limits_{k=0}^{m}H_k\triangle^k u(x)-u(y)}{|x-y|^{n+2t}}dy\\
&=& C_{n,t}\lim_{\varepsilon \ra 0}\int_{R^n\backslash B_\varepsilon(x)}
\frac{\mathop\Sigma\limits_{k=0}^{m}H_k\triangle^k u(x)-u(y)}{|x-y|^{n+2t}}dy,
\end{eqnarray*}
where
$$H_k=\frac{1}{2^k k!n(n+2)\cdots(n+2k-2)}.$$
An equivalent approach to define it is via the  difference quotient:
$$
(-\lap)^{t}u(x)=
C_{n,t}\int_{R^n}\frac{\mathop\Sigma\limits_{k=0}^{2m+2}
(-1)^{k-m-1}C_{2m+2}^k u(x+(k-m-1)y)}{|y|^{n+2t}}dy,
$$
with
$$C_{2m+2}^k=\frac{(2m+2)!}{k!(2m+2-k)!}.$$

From the defining integral in (\ref{2016672}), we
 can see that both $(-\lap)^t$ and
 $(-\lap)^{\alpha/2}$ are nonlocal pseudo-differential operators.
Such non-locality brings in the real difference between the usual Laplacians and the fractional Laplacians, and poses a strong barrier
in the generalization of many useful results from the Laplacians to
the fractional ones. For example, when considering a Dirichlet
problem
 involving $-\lap$ in a bounded domain $\Omega$,
    it is
 sufficient to know the behavior
of solution on $\partial\Omega$. But when it comes to
 $(-\lap)^{\alpha/2}$, we need to gather information
  on both $\partial\Omega$ and
$R^n\backslash\Omega$.

Due to its non-local nature, the fractional Laplacian have
been receiving intense attention from researchers in physics,
astrophysics, mechanics, mathematics and economics.
Scientists have been using integrations and differentiation of fractional orders to describe the  behavior of objects and systems.
 In physics, it is used to derive heat kernel estimates for a large class of symmetric jump-type processes (see \cite{CKK},
 \cite{BBCK}). Its application is also seen in the study of
 the acoustic wave equation, which was an important point of reference in the development of the electromagnetic wave equation. In astrophysics, researchers use it to model the dynamics in the Hamiltonian chaos (see \cite{Z}).
In probability, it is defined
as the generator of $\alpha$-stable L$\acute{e}$vy processes
that represents random motions, such as the Brownian motion and the
Poisson process(see \cite{D} and \cite{K}). In finance,
 it models jump precess (see \cite{RP}). For more details,
 please see \cite{B}, \cite{BGR}, \cite{BG} and the references there in.

In this paper, we consider the Navier boundary value problem:
\begin{equation}
 \left\{\begin{array}{ll}
(-\lap)^t u(x)=u^p(x),& x\in R^n_+,\\
(-\lap)^i u(x)=0,& x\notin R^n_+,\\
\end{array} \right.
\label{20165111}
\end{equation}
with $i=0, 1, 2, \cdot\cdot\cdot, m$.
We want to derive the nonexistence of positive solutions.

The nonexistence result is often famously known as the Liouville-type theorem.
Interestingly enough, in PDE analysis,
by applying the nonexistence to blowed-up-and-rescaled equations, one is able to derive
the a priori estimates for the solutions of the original equations. Then,  combining the estimates with methods like the topological degree and
continuation method, it forms a powerful tool to derive the existence of solutions.

Considerable work have been devoted to the
study of the Liouville theorems associated with the
usual fractional Laplacian $(-\lap)^{\alpha/2}$ and interesting results have been obtained.

In \cite{ZCCY}, the authors proved that the only solution for
$$ \left\{\begin{array}{ll}
(-\lap)^{\alpha/2} u(x) = 0 , \;\; x \in \mathbb{R}^n ,\\
u(x)\geq 0,  \;\; x \in \mathbb{R}^n,
\end{array}
\right.
$$
is constant.

In \cite{CFY}, Chen, Fang and Yang studied
the following Dirichlet problem involving the fractional Laplacian:
\begin{equation*}
\left\{\begin{array}{ll}
(- \lap)^{\alpha/2} u = u^p, & x \in \mathbb{R}^n_+, \\
u \equiv 0 , & x \not{\in} \mathbb{R}^n_+ .
\end{array}
\right.
\end{equation*}
By studying its equivalent integral equation
$$
 u(x) = \int_{\mathbb{R}^n_+} G(x, y) u^p(y) d y,
$$
they obtained the non-existence of positive solutions in the critical and subcritical cases $1 < p \leq \frac{n+\alpha}{n-\alpha}$
 under no restrictions on the growth of the solutions.

In \cite{CDY}, the authors considered an
$\alpha$-harmonic problem:
$$ \left\{\begin{array}{ll}
(-\lap)^{\alpha/2} u(x) = 0 , \;\; x \in \mathbb{R}^n ,\\
\displaystyle\underset{|x| \ra \infty}{\underline{\lim}} \frac{u(x)}{|x|^{\gamma}} \geq 0,
\end{array}
\right.
$$
for some $0 \leq \gamma \leq 1$ and $\gamma < \alpha$.
There they proved that
 $u$ must be constant throughout $\mathbb{R}^n$.

For more results on the Liouville theorems, please see \cite{Fa}, \cite{L}, \cite{ZLL} and the references therein). To our best knowledge, the Liouville theorems concerning
the higher-order fractional Laplacian are few and new.
In  \cite{LZ},  among which, we
derived the equivalent integral form of the above differential equation:

\begin{mlem}\label{2016662}
Suppose $u$ is a positive solution of (\ref{20165111}).
If $(-\Delta)^m u\geq 0$ and $2m<n$, then u must satisfy the integral equation
\begin{equation}
u(x)=\int_{R^n_+}G_{2t}^+(x,y)u^p(y)dy.\label{0.4}
\end{equation}
And vice versa.
\end{mlem}
Based on the equivalence result, at the end of the paper we conjectured that

\emph{Assume that  $u \in L^{\frac{n(p-1)}{2t}}_{loc}(R^n_+)$ for $p>1$ is
a nonnegative solution to  (\ref{20165111}). Then $u$ must be trivial.}

As a continuation of \cite{LZ}, here we
partially prove the conjecture under a stronger assumption that
 $u \in L^{\frac{n(p-1)}{2t}}(R^n_+)$.
Let $R^+=\{x=(x_1,\dots,x_n)\in R^n \mid x_n>0\}$.
\begin{mthm}\label{201651616}
Assume that $(-\Delta)^m u\geq 0$ and $2m<n$. If $u\in C_{loc}^{2m,1} \cap L_{2t}$ is a nonnegative solution for
\begin{equation}
u(x)=\int_{R^n_+}G_{2t}^+(x,y)u^p(y)dy,
\end{equation}
then $u$ must be trivial.
\end{mthm}

As an immediate consequence of Lemma \ref{2016662} and Theorem
\ref{201651616}, we have

\begin{mthm}
  For
$p>1$ and $n>2t$, suppose
 $u \in L^{\frac{n(p-1)}{2t}}(R^n_+)$ is a nonnegative solution for
 (\ref{20165111}). If $(-\lap)^m u\geq 0$ in $R^n_+$, $u$ must be trivial.
\end{mthm}

The pseudo-differential operator $(-\lap)^{t}$
can also be defined inductively (see \cite{OS}) as:
\be\label{2016672}
(-\lap)^{t}=(-\lap)^{\alpha/2}\circ (-\Delta)^m,
\ee
where $(-\lap)^{\alpha/2}$ is the fractional Laplacian.
This allows us to split a single higher-order fractional equation into a system. We will develop an iteration method
based on the following narrow region principle to deal with
the differential problem directly.

\begin{mthm}
Assume that $\Omega$ is a bounded narrow region in
 $\Sigma_{\la}=\{x\in \mathbb{R}^n  \mid x_1<\lambda\}$.
Without loss of generality, we may assume that $\Omega$ is contained in the slab $\{x\in \mathbb{R}^n \mid \la-l<x_1<\la\}$ with $\varepsilon>0$ small.
Consider
\begin{equation*}
\left\{\begin{array}{ll}
-\lap U(x)+c_1(x)V(x)\geq 0, &x \in \Omega,\\
(-\lap)^{\alpha/2} V(x)+c_2(x) U(x)\geq0,&x \in \Omega,\\
 V(x^{\la})=-V(x), & x \in \Sigma_{\la}, \\
  U(x^{\la}a)=-U(x), & x \in \Sigma_{\la}, \\
U(x),\,V(x)\geq0,& x \in \Sigma_{\la}\backslash \Omega,
\end{array}
\right.
\end{equation*}
where $c_i(x)\leq 0$ in $\Omega$  and are bounded for
$i=1,\: 2$, $U \in C^2$ and $V \in C_{loc}^{1,1}(\Omega \cap L_\alpha)$ are lower semi-continuous in $\bar{\Omega}$.

Then for $l$ sufficiently small, we have
$$
U(x), \:V(x)\geq0,\; x \in \Omega.
$$
For unbounded $\Omega$, the above result still holds  on condition that
$$U(x), \, V(x) \ra 0, \quad |x| \ra \infty.$$
Further, if either $U(x)$ or $V(x)$ equals 0 at some point in $\Omega$, then
$$
U(x),\,V(x)\equiv 0,\; x \in R^n.
$$
\end{mthm}

Then  we are able to prove that
\begin{mthm}\label{2016913}
Let $t=1+\alpha/2$ for $0<\alpha<2$. Assume that $u \in C_{loc}^{3,1}(B_1(0))$
is a positive solution of
\begin{equation*}
\left\{\begin{array}{ll}
(-\lap)^t u=f(u),  &x \in B_1,\\
u=\lap u=0, &x \in R^{n}\backslash B_1,
\end{array}
\right.
\end{equation*}
where $f(t)$ is Lipschitz continuous and increasing in $t$. Then $u$ must be symmetric about the origin.
\end{mthm}

With the symmetry result, one can
continue to carry out estimates on the
regularity results and derive useful Sobolev inequalities.
For readers who are interested please see \cite{CLO1}
, \cite{CLO2} and the references therein.

This paper is organized as follows: In Section 2, we obtain some
essential inequalities for the Green's function related to
$(-\Delta)^t$; In Section 3, we prove Theorem
\ref{201651616} through the method of moving planes in
integral forms. We close the paper with a proof of the narrow region principle and Theorem
\ref{2016913}.

\section{Properties of the Green's Function}

\smallskip
First we derive some properties of the Green's function
$G_{2t}^+(x,y)$, which is essential in the process of moving the planes.

For any real number $\lambda$, let
$$\Sigma_{\lambda}=\left\{x=(x',x_n)\in R^n_+\mid 0<x_n<\lambda \right\},$$
$$ T_{\lambda}=\left\{x\in
R^n_+\mid x_n=\lambda \right\},$$ and let
$$x^{\lambda}=(x_1,x_2,\cdots,2\lambda-x_n)$$
be the reflection of the point $x=(x',x_n)$ abo ut the
plane $T_{\lambda}$.

\begin{lem}\label{20165113}
The Green's function $G_{2t}^+(x,y)$ satisfies the following properties:
\begin{enumerate}
  \item For any $x,y\in \Sigma_{\lambda}$ and $x\neq y$, we have
\begin{equation}\label{20165161}
G_{2t}^+(x^{\lambda},y^{\lambda})>\max\{G_{2t}^+(x^{\lambda},y),
G_{2t}^+(x,y^{\lambda})\}
\end{equation}
and
\begin{equation}\label{20165162}
G_{2t}^+(x^{\lambda},y^{\lambda})-G_{2t}^+(x,y)
>|G_{2t}^+(x^{\lambda},y)-G_{2t}^+(x,y^{\lambda})|.
\end{equation}

  \item For any $x\in \Sigma_{\lambda}$, $y\in \Sigma_{\lambda}^{C}$,
   it holds
\begin{equation}\label{20165163}
G_{2t}^+(x^{\lambda},y)> G_{2t}^+(x,y) \; \mbox{ and }
G_{2t}^+(y, x^{\lambda}) > G_{2t}^+(y,x).
\end{equation}
  \item For $ x, \; y \in R^n_+ $ ,
  \begin{equation}\label{20165165}
  G_{2t}^+(x,y)\leq \frac{C}{|x-y|^{n-2t}}.
\end{equation}
\end{enumerate}
\end{lem}

\textbf{Proof.}\; It's well known that the above inequalities are
true for the Green's functions associated with $-\lap$ and $(-\lap)^m$. In the proof below, we will make use of this fact.

Let $ \tilde{\Sigma}_\la$ be the reflection of $\Sigma_\la$ about $T_\la$, and
$R^n_+=\Sigma_\la \cup \tilde{\Sigma}_\la \cup D_\la$.
To prove the lemma, we recall that
 $G_{2m}^+(x,y)$ and $ G_{\alpha}^+(x,y)$
satisfy inequality (\ref{20165161})-(\ref{20165163}) as well. Therefore,
\begin{enumerate}
  \item  \begin{eqnarray}
&&G_{2t}^+(x^{\lambda},y^{\lambda})-G_{2t}^+(x^{\lambda},y) \nonumber\\
 &=& \int_{R^n_+} G_{2m}^+(x^{\lambda}, z)
            \big(G_{\alpha}^+(z, y^{\lambda}) -  G_{\alpha}^+(z, y)\big) d z  \nonumber\\
&=& \int_{\Sigma_\la \cup \tilde{\Sigma}_\la \cup D_\la} G_{2m}^+(x^{\lambda}, z)
            \big(G_{\alpha}^+(z, y^{\lambda}) -  G_{\alpha}^+(z, y)\big) d z   \nonumber\\
&=& \int_{\Sigma_{\lambda}}\bigg(G_{2m}^+(x^{\lambda}, z)
            \big(G_{\alpha}^+(z, y^{\lambda}) -  G_{\alpha}^+(z, y)\big) \nonumber\\
            &&+ G_{2m}^+(x^{\lambda}, z^{\lambda})
  \big(G_{\alpha}^+(z^{\lambda}, y^{\lambda}) -  G_{\alpha}^+(z^{\lambda}, y)\big) \bigg)d z  \nonumber\\
&&+ \int_{D_{\lambda}} G_{2m}^+(x^{\lambda}, z)
            \big(G_{\alpha}^+(z, y^{\lambda}) -  G_{\alpha}^+(z, y)\big) d z  \nonumber\\
&=& I_1 + I_2 .
 \label{20165164}
 \end{eqnarray}
Since $G_{\alpha}^+$ satisfies (\ref{20165162}), we have
$$G_{\alpha}^+(z^{\lambda}, y) -  G_{\alpha}^+(z, y)\geq
G_{\alpha}^+(z, y^{\lambda}) -  G_{\alpha}^+(z^{\lambda}, y^{\lambda}).$$
Thus
\begin{eqnarray*}
I_1&=& \int_{\Sigma_{\lambda}}  \bigg(G_{2m}^+(x^{\lambda}, z)
            \big(G_{\alpha}^+(z, y^{\lambda}) -  G_{\alpha}^+(z, y)\big)\\
 &&+  G_{2m}^+(x^{\lambda}, z^{\lambda})
  \big(G_{\alpha}^+(z^{\lambda}, y^{\lambda}) -  G_{\alpha}^+(z^{\lambda}, y)\big)\bigg) d z  \\
   &\geq & \int_{\Sigma_{\lambda}}  \bigg(- G_{2m}^+(x^{\lambda}, z)[ G_{\alpha}^+(z^{\lambda}, y^{\lambda}) - G_{\alpha}^+(z^{\lambda}, y)] \\
   &&+   G_{2m}^+(x^{\lambda}, z^{\lambda})[ G_{\alpha}^+(z^{\lambda}, y^{\lambda}) - G_{\alpha}^+(z^{\lambda}, y)]  \bigg)d z  \nonumber\\
&=& \int_{\Sigma_{\lambda}} [  G_{2m}^+(x^{\lambda}, z^{\lambda}) - G_{2m}^+(x^{\lambda}, z)][ G_{\alpha}^+(z^{\lambda}, y^{\lambda}) - G_{\alpha}^+(z^{\lambda}, y)]  d z  \nonumber\\
&> & 0. \label{20165166}
\end{eqnarray*}
Because (\ref{20165162}) is true for $G_{\alpha}^+$, it implies that
 $$ I_2=\int_{D_{\lambda}} G_{2m}^+(x^{\lambda}, z)
\big(G_{\alpha}^+(z, y^{\lambda}) -  G_{\alpha}^+(z, y)\big) d z  >0.$$
Together with (\ref{20165166}), it yields
$$G_{2t}^+(x^{\lambda},y^{\lambda})\geq G_{2t}^+(x^{\lambda},y).$$
  Similarly, one can show that
 $$G_{2t}^+(x^{\lambda},y^{\lambda})\geq G_{2t}^+(x,y^{\lambda}).$$
 This proves (\ref{20165161}). To prove (\ref{20165162}),
 let
 \begin{eqnarray}
& &[ G_{2t}^+(x^{\lambda}, y^{\lambda}) - G_{2t}^+(x,y) ] -
 [ G_{2t}^+(x^{\lambda},y) - G_{2t}^+(x, y^{\lambda})] \nonumber \\
&=& \int_{\Sigma_{\lambda}} K dz +   \int_{D_{\lambda}} Q dz , \label{20165167}
\end{eqnarray}
  where
\begin{eqnarray}
K &=& [G_{2m}^+(x^{\lambda},z) + G_{2m}^+(x,z)]
[G_{\alpha}^+(z,y^{\lambda})-G_{\alpha}^+(z,y)] \nonumber \\
\smallskip
&&+[G_{2m}^+(x^{\lambda},z^{\lambda})+ G_{2m}^+(x,z^{\lambda})]
[G_{\alpha}^+(z^{\lambda},y^{\lambda})-G_{\alpha}^+(z^{\lambda},y)] \nonumber \\
\smallskip
&>& -[G_{2m}^+(x^{\lambda},z) + G_{2m}^+(x,z)]
[G_{\alpha}^+(z^{\lambda},y^{\lambda})-G_{\alpha}^+(z^{\lambda},y)] \nonumber \\
\smallskip
&&+[G_{2m}^+(x^{\lambda},z^{\lambda})+ G_{2m}^+(x,z^{\lambda})]
[G_{\alpha}^+(z^{\lambda},y^{\lambda})-G_{\alpha}^+(z^{\lambda},y)] \nonumber \\
\smallskip
&=& [G_{2m}^+(x^{\lambda},z^{\lambda})+ G_{2m}^+(x,z^{\lambda})
-G_{2m}^+(x^{\lambda},z) - G_{2m}^+(x,z)]\nonumber \\
\smallskip
&&\cdot [G_{\alpha}^+(z^{\lambda},y^{\lambda})-G_{\alpha}^+(z^{\lambda},y)] \nonumber \\
& > & 0 . \label{20165168}
\end{eqnarray}
By (\ref{20165163}), we also have
\begin{equation}
Q = [G_{2m}^+(x^{\lambda},z)+G_{2m}^+(x,z)]
[G_{\alpha}^+(z,y^{\lambda})-G_{\alpha}^+(z,y)] > 0 .
\label{20165169}
\end{equation}

Combining (\ref{20165167}), (\ref{20165168}), and (\ref{20165169}), we arrive at
$$ G_{2t}^+(x^{\lambda}, y^{\lambda}) - G_{2t}^+(x,y)
>  G_{2t}^+(x^{\lambda},y) - G_{2t}^+(x, y^{\lambda}) .$$
Similarly, one can show that
$$ G_{2t}^+(x^{\lambda}, y^{\lambda}) - G_{2t}^+(x,y)
>  G_{2t}^+(x,y^{\lambda}) - G_{2t}^+(x^{\lambda}, y).$$
This verifies (\ref{20165162}).

  \item For $x \in \Sigma_{\lambda}$ and $y \in \Sigma_{\lambda}^C $, let
\begin{equation}
G_{2t}^+(x^{\lambda},y)-G_{2t}^+(x,y) = \int_{\Sigma_{\lambda}} M dz +   \int_{D_{\lambda}} N dz ,
\label{201651610}
\end{equation}
where
\begin{eqnarray}
M &=& [G_{2m}^+(x^{\lambda}, z)-G_{2m}^+(x,z)]G_{\alpha}^+(z,y)\nonumber \\
&&+ [G_{2m}^+(x^\lambda, z^\lambda)-G_{2m}^+(x,z^\lambda)]G_{\alpha}^+(z^\lambda,y) \nonumber \\
&>& [G_{2m}^+(x^\lambda, z^\lambda)-G_{2m}^+(x,z^\lambda)]
[G_{\alpha}^+(z^\lambda,y)-G_{\alpha}^+(z,y)] > 0 . \label{201651611}
\end{eqnarray}
Meanwhile, we have
\begin{equation}
N = [G_{2m}^+(x^\lambda,z)-G_{2m}^+(x,z)]G_{\alpha}^+(z,y) > 0 .
\label{201651612}
\end{equation}
Combining (\ref{201651610}), (\ref{201651611}) and (\ref{201651612}),
we arrive at (\ref{20165163}).

  \item It's well known that for $0<\alpha\leq2$, $x,\; y\in \mathbb{R}^n_+$,
$$
G_\alpha^+(x,y)=\frac{A_{n,\alpha}}{s^{(n-\alpha)/2}}\left[1- \frac{B_{n,\alpha}}{(s+t)^{(n-2)/2}} \int_0^{\frac{s}{t}}
\frac{(s-tb)^{(n-2)/2}}{b^{\alpha/2}(1+b)} d b \right],
$$
is the Green's function for $(-\lap)^{\alpha/2}$ in $\mathbb{R}^n_+$.
 Here
$$s = |x-y|^2 \;\; \mbox{ while } \;\; t = 4 x_n y_n .$$
We know that
$$G_{2m}^+(x,y)=\frac{1}{|x-y|^{n-2m}}-\frac{1}{|x^*-y|^{n-2m}},$$
thus
\begin{equation*}
  G_{2m}^+(x,y)\leq \frac{C}{|x-y|^{n-2m}}, \;\; x, y \in R^n_+,
\end{equation*}
and
\begin{equation*}
   G_\alpha^+(x,y)\leq \frac{C}{|x-y|^{n-\alpha}}, \;\; x, y \in R^n_+.
\end{equation*}
Hence
\begin{eqnarray*}
 &&G_{2t}^+(x,y) \\
 &=& \int_{R^n_+}G_{2m}^+(x,z)G_{\alpha}^+(z,y)dz \\
   &\leq&  \int_{R^n_+} \frac{C}{|x-z|^{n-2m}|z-y|^{n-\alpha}}dz\\
   &=& \int_{\bar{z}_n> -y_n} \frac{C}{|x-y-\bar{z}|^{n-2m}|\bar{z}|^{n-\alpha}},
  \;\;\bar{z}+y=z\\
   &=& \frac{C}{|x-y|^{n-2t}}\int_{\tilde{z}_n>-\frac{y_n}{|x-y|}}
    \frac{d\tilde{z}}{|\tilde{z}|^{n-\alpha}|\frac{x-y}{|x-y|}-\tilde{z}|^{n-2m}}d \bar{z},\;\;\bar{z}=|x-y|\tilde{z}\\
   &=& \frac{C}{|x-y|^{n-2t}}.
\end{eqnarray*}
\end{enumerate}

\smallskip

\section{Nonexistence of Positive Solutions for the Integral Equation}

\textbf{Proof of Theorem \ref{201651616}}
\smallskip

Let $$u_\la(x)=u(x^\la),\quad w_\la(x)=u_\la(x)-u(x).$$

Then
\begin{lem}
\be\label{20165114}
w_\la(x)\geq \int_{\Sigma_\la}[G_{2t}(x^\la, y^\la)- G_{2t}(x, y^\la)]
(u_\la^p(y)-u^p(y))\,dy.
\ee
\end{lem}

\textbf{Proof}
By (\ref{0.4}) and Property 2 in Lemma \ref{20165113}, we have
\begin{eqnarray*}\nonumber
w_\la(x)&=& u_\la(x)-u(x) \\\nonumber
  &=& \int_{R^n_+}G_{2t}^+(x^\la,y)u^p(y)dy- \int_{R^n_+}G_{2t}^+(x,y)u^p(y)dy
  \\\nonumber
 &=& \int_{\Sigma_\la \cup \tilde{\Sigma}_\la \cup D_\la}G_{2t}^+(x^\la,y)u^p(y)dy- \int_{\Sigma_\la \cup \tilde{\Sigma}_\la \cup D_\la}G_{2t}^+(x,y)u^p(y)dy
  \\\nonumber
 &=&\int_{\Sigma_\la}[G_{2t}^+(x^\la,y)-G_{2t}^+(x,y)]u^p(y)dy\\
    &&+\int_{\Sigma_\la}[G_{2t}^+(x^\la,y^\la)-G_{2t}^+(x,y^\la)]u^p_\la(y)dy \\
 &&+\int_{D_\la}[G_{2t}^+(x^\la,y)- G_{2t}^+(x,y) ]u^p(y)dy\\
   &\geq & \int_{\Sigma_\la}[G_{2t}^+(x^\la,y)-G_{2t}^+(x,y)]u^p(y)dy\\
    &&+\int_{\Sigma_\la}[G_{2t}^+(x^\la,y^\la)-G_{2t}^+(x,y^\la)]u^p_\la(y)dy.
    \end{eqnarray*}
 By Property 1 in Lemma \ref{20165113}, one can see that for $x, \:y \in \Sigma_\la$,
 $$G_{2t}^+(x^\la,y)-G_{2t}^+(x,y)\geq
 G_{2t}^+(x^\la,y^\la)-G_{2t}^+(x,y^\la).$$
 Therefore,
  \begin{eqnarray*}
  w_\la(x) &\geq&  \int_{\Sigma_\la}[G_{2t}^+(x,y^\la)-G_{2t}^+(x^\la,y^\la)]u^p(y)dy\\
  &&     +\int_{\Sigma_\la}[G_{2t}^+(x^\la,y^\la)-G_{2t}^+(x,y^\la)]u^p_\la(y)dy \\
   &\geq& \int_{\Sigma_\la} [G_{2t}^+(x^\la,y^\la)-G_{2t}^+(x,y^\la)]
   [u^p_\la(y)-u^p(y)]dy.
\end{eqnarray*}
This proves the lemma.

Next we use the method of moving planes to derive a contradiction assuming that
(\ref{20165111}) has positive solutions.

\emph{Step 1. Start moving the plane $T_\la$ from near $x_n=0$ to the right along
the $x_n$ axis.}

Let $\Sigma_\la^-=\{x \in \Sigma_\la \mid w_\la(x)<0   \}$.
We show that
\be\label{20165115}
w_\la(x)\geq 0, \quad a.e. \; \Sigma_\la.
\ee

If not, then for any $x^o \in \Sigma_\la^-$, by (\ref{20165114}) and the
\textit{mean value theorem}, one have
\begin{eqnarray*}
0<-w_\la(x^o) &\leq &
\int_{\Sigma_\la^-}[G_{2t}(x^\la, y^\la)- G_{2t}(x, y^\la)]
(u_\la^p(y)-u^p(y))\,dy\\
   &\leq &C \int_{\Sigma_\la^-}  \frac{1}{|x^\la-y^\la|^{n-2t}}
   (u_\la^p(y)-u^p(y))\,dy\\
   &\leq &C \int_{\Sigma_\la^-}  \frac{1}{|x-y|^{n-2t}}
   p u^{p-1}(y)(-w_\la(y))\,dy.
\end{eqnarray*}

Next we need an equivalent form of the Hardy-Little-Sobolev inequality.
\begin{lem}
 Assume $0<\alpha<n$ and $\Omega\subset \mathbb{R}^n$. Let $g\in L^{\frac{np}{n+\alpha p}}(\Omega)$ for $\frac{n}{n-\alpha}<p<\infty.$ Define
$$
Tg(x):=\int_{\Omega}\frac{1}{|x-y|^{n-\alpha}}g(y)dy.
$$
Then
\begin{equation}\label{n8}
 \|Tg\|_{L^{p}(\Omega)}\leq C(n,p, \alpha)\|g\|_{L^{\frac{np}{n+\alpha p}}(\Omega)}.
\end{equation}
\label{lemhls}
 \end{lem}
The proof of this lemma is standard and can be found in \cite{CL} or \cite{CL1}.

By Lemma \ref{lemhls} and the H$\ddot{o}$lder inequality, it follows that
\begin{eqnarray*}
\|w_\la\|_{L^q(\Sigma_\la^-)} &\leq &
 C\| u^{p-1}w_\la \|_{L^{\frac{nq}{n+2tq}}(\Sigma_\la^-)}\\
  &\leq& C\| u^{p-1}\|_{L^\frac{n}{2t}(\Sigma_\la^-)}  \|w_\la\|_{L^q(\Sigma_\la^-)},
\end{eqnarray*}
for any $\frac{n}{n-2t}<q<\infty$.
Since $u \in  L^{\frac{n(p-1)}{2t}}(R^n_+)$, for $\la$  sufficiently negative,
 we have
\begin{equation}\label{20165118}
C\| u^{p-1}\|_{L^\frac{n}{2t}(\Sigma_\la^-)}<\frac{1}{2}.
\end{equation}
Thus
$$\|w_\la\|_{L^q(\Sigma_\la^-)} < \frac{1}{2} \|w_\la\|_{L^q(\Sigma_\la^-)},$$
i.e. $\Sigma_\la^-$ is measure 0.

This proves (\ref{20165115}) and completes step 1.

\emph{Step 2. Keep moving the plane until the limiting position}
$$\la_0=\sup \{ \lambda\leq\infty  \mid  w_{\mu}(x)\geq
0,\,\,\,\forall x\in\Sigma_{\mu},\,\,\,\mu\leq\lambda   \}.$$

We claim that
\be\label{20165131}
\la_0=\infty.
\ee
If not, then for $\la_0<\infty$, we can show that
\be\label{20165132}
w_{\la_0}(x)\equiv0.
\ee
Thus for any $x^0 \in \partial R^n_+$, we have
$$0<u_{\la_0}(x^0)=u(x^0)=0.$$
The contradiction establishes (\ref{20165131}).

To prove (\ref{20165132}), we suppose its contrary is true, i.e.
there exists a non-measure-zero
set
$D \subseteq \Sigma_{\la_0}$ such that
\be
w_{\la_0}(x)>0, \quad x \in D.
\ee
Thus, for any $x \in \Sigma_{\la_0}$, by (\ref{20165114}),
\begin{eqnarray*}
w_{\la_0}(x) &\geq & \int_{\Sigma_{\la_0}}[G_{2t}(x^{\la_0}, y^{\la_0})- G_{2t}(x, y^{\la_0})]
(u_{\la_0}^p(y)-u^p(y))\,dy \\
   &\geq &  \int_D [G_{2t}(x^{\la_0}, y^{\la_0})- G_{2t}(x, y^{\la_0})]
(u_{\la_0}^p(y)-u^p(y))\,dy \\
   &>& 0.
\end{eqnarray*}
Such strict positivity allows us to keep moving $T_{\la_0}$ to the right while
preserving (\ref{20165115}). In other words, there exists some $\varepsilon>0$ small
such that for any $\la \in (\la_0, \la_0+\varepsilon)$, it holds that
\be\label{20165117}
 w_{\la}(x)\geq 0, \quad a.e.\; \Sigma_\la.
 \ee
This is a contradiction with the definition of $\la_0$.
It thus verifies (\ref{20165132}).

To prove (\ref{20165117}), for any small $\eta>0$, we can choose $R$ sufficiently large, so
that
\begin{eqnarray}\left(\int_{R^n\setminus B_R(0)}u^{\tau +1}(y)dy\right)^{\frac{\alpha}{n}}<\eta .
\label{6}\end{eqnarray}

Fix this $R$, we will show that the measure of $\Sigma_{\lambda}^- \cap B_R(0)$ is sufficiently small
as $\lambda$ close to $\lambda_o$.

By the definition of $\la_o$ and Lemma \ref{20165114}, it is
trivial to deduce that $w_{\lambda_o} (x)> 0$ in the
interior of $\Sigma_{\lambda_o}$. For any $\delta > 0$, let
$$ E_{\delta} = \{ x \in \Sigma_{\lambda_o} \cap B_R(0) \mid w_{\lambda_o} (x) > \delta \} \;\; \mbox{ and } \; \;
F_{\delta} = \left(\Sigma_{\lambda_o} \cap B_R(0)\right) \setminus E_{\delta} .$$ Then obviously
$$ \lim_{\delta \ra 0} \mu (F_{\delta}) = 0 .$$

For $\lambda > \lambda_o$, let
$$ D_{\lambda} = \left( \Sigma_{\lambda} \setminus \Sigma_{\lambda_o} \right) \cap B_R(0) .$$
Then it is easy to see that
\begin{equation}
\left(\Sigma_{\lambda}^- \cap B_R(0)\right) \subset \left(\Sigma_{\lambda}^- \cap E_{\delta}\right) \cup F_{\delta} \cup D_{\lambda} .
\label{AB1}
\end{equation}
Apparently, the measure of $D_{\lambda}$ gets small as $\lambda$ approaches $\lambda_o$. We show that the measure of $\Sigma_{\lambda}^- \cap E_{\delta}$ can also be sufficiently small as $\lambda$ close to $\lambda_o$. In fact, for any $x \in \Sigma_{\lambda}^- \cap E_{\delta}$, we have
$$ w_{\lambda} (x) = u(x^{\lambda}) - u(x^{\lambda_o}) + u(x^{\lambda_o}) - u(x) < 0 .$$
Hence
$$ u(x^{\lambda_o}) - u(x^{\lambda}) > w_{\lambda_o}(x) > \delta .$$
It follows that
\begin{equation}
\left(\Sigma_{\lambda}^- \cap E_{\delta}\right) \subset G_{\delta} \equiv \{ x \in B_R(0) \mid u(x^{\lambda_o}) - u(x^{\lambda}) > \delta \} .
\label{AB2}
\end{equation}
While by the well-known Chebyshev inequality, we have
\begin{eqnarray*}
\mu (G_{\delta}) & \leq & \frac{1}{\delta^{\tau +1}} \int_{G_{\delta}} |u(x^{\lambda_o}) - u(x^{\lambda})|^{\tau +1} d x  \\
 &\leq&\frac{1}{\delta^{\tau +1}} \int_{B_R(0)} |u(x^{\lambda_o}) - u(x^{\lambda})|^{\tau +1} d x .
\end{eqnarray*}
For each fixed $\delta$, as $\lambda$ close to $\lambda_o$, the right hand side of the above inequality can be made as small as we wish. Therefore by (\ref{AB2}) and (\ref{AB1}), the measure of $\Sigma_{\lambda}^- \cap B_R(0)$ can also be sufficiently small. Combining this with
(\ref{6}), we arrive at (\ref{20165117}).

It follows from  $\la_0=\infty$ that $u(x)$ is monotone increasing in the
$x_n$ direction. This violates our assumption that
$u \in L^{\frac{n(p-1)}{2t}}(R^n_+)$. Therefore, (\ref{20165111})
has no positive solution.

This completes the proof of the theorem.

\section{Symmetry of Positive Solutions}

To derive the symmetry through the iteration method, we first prove a key ingredient -- the narrow region principle.

\subsection{Narrow Region Principle}

\begin{thm}
Assume that $\Omega$ is a bounded narrow region in
 $\Sigma_{\la}=\{x\in \mathbb{R}^n  \mid x_1<\lambda\}$.
Without loss of generality, we may assume that $\Omega$ is contained in the slab $\{x\in \mathbb{R}^n \mid \la-l<x_1<\la\}$ with $\varepsilon>0$ small.

Consider
\begin{equation}\label{2016541}
\left\{\begin{array}{ll}
-\lap U(x)+c_1(x)V(x)\geq 0, &x \in \Omega,\\
(-\lap)^{\alpha/2} V(x)+c_2(x) U(x)\geq0,&x \in \Omega,\\
 V(x^{\la})=-V(x), & x \in \Sigma_{\la}, \\
  U(x^{\la}a)=-U(x), & x \in \Sigma_{\la}, \\
U(x),\,V(x)\geq0,& x \in \Sigma_{\la}\backslash \Omega,
\end{array}
\right.
\end{equation}
where $c_i(x)\leq 0$ in $\Omega$  and are bounded for
$i=1,\: 2$, $U \in C^2$ and $V \in C_{loc}^{1,1}(\Omega \cap L_\alpha)$ are lower semi-continuous in $\bar{\Omega}$.

Then for $l$ sufficiently small, we have
\be\label{20165126}
U(x), \:V(x)\geq0,\; x \in \Omega.
\ee

For unbounded $\Omega$, (\ref{20165126}) still holds  on condition that
$$U(x), \, V(x) \ra 0, \quad |x| \ra \infty.$$

Further, if either $U(x)$ or $V(x)$ equals 0 at some point in $\Omega$, then
\be\label{20165124}
U(x),\,V(x)\equiv 0,\; x \in R^n.
\ee
\end{thm}

\textbf{Proof. } 
If $V(x)\geq 0$ is not true, then there must exist some
$ x^0\in \Omega$ such that
\begin{equation*}
V(x^0)=\min_{\Sigma_{\la}}V(x) <0.
\end{equation*}
Then
\begin{eqnarray}
(-\lap)^{\alpha/2}V(x^0)&=&C_{n,\alpha}PV\int_{\mathbb{R}^n}
\frac{V(x^0)-V(y)}
{|x^0-y|^{n+\alpha}}dy\nonumber\\
&=&C_{n,\alpha}PV\left\{\int_{\Sigma_{\la}} \frac{V(x^0)-V(y)}
{|x^0-y|^{n+\alpha}}dy+\int_{R^n\backslash {\Sigma_{\la}}} \frac{V(x^0)-V(y)}
{|x^0-y|^{n+\alpha}}dy \right\}\nonumber\\
&=&C_{n,\alpha}PV\left\{\int_{\Sigma_{\la}} \frac{V(x^0)-V(y)}
{|x^0-y|^{n+\alpha}}dy+\int_{\Sigma_{\la}} \frac{V(x^0)+V(y)}
{|x^0-y^{\la}|^{n+\alpha}}dy\right\}\nonumber\\
&\leq&
C_{n,\alpha}\int_{\Sigma_{\la}}
\left\{\frac{V(x^0)-V(y)}{|x^0-y^{\la}|^{n+\alpha}}
+\frac{V(x^0)+V(y)}{|x^0-y^{\la}|^{n+\alpha}}\right\}dy\nonumber\\
&=&C_{n,\alpha}\int_{\Sigma_{\la}} \frac{2V(x^0)}
{|x^0-y^{\la}|^{n+\alpha}}dy\nonumber\\
&\leq&CV(x^0)\int_{B_1(x^o)\setminus B_l(x^o)}\frac{1}{|x^o-y|^{n+\alpha}}dy\nonumber\\
&\leq &\frac{CV(x^0)}{l^{\alpha}}.\label{2016527}
\end{eqnarray}
Together with
$$(-\lap)^{\alpha/2} V(x)+c_2(x) U(x)\geq0,$$
we deduce that $$U(x^0)<0.$$
It thus implies that there exists some
 $\bar{x} \in \Sigma_{\la} \cap B_1(0)$ such that
\begin{equation}\label{20165121}
U(\bar{x})=\min_{\Sigma_{\la}}U(x)<0.
\end{equation}

For $\delta>0$ small,  let
$ \phi(x)=\sin\bigg( \frac{x_1-\la+l+\delta}{l} \bigg)$.
Then $ \phi(x)$ has positive bounds and satisfies
$\lap \phi(x)= -\frac{\phi(x)}{l^2}$.
Let $w(x)=\frac{U(x)}{\phi(x)}$. It follows from (\ref{20165121}) that
there exists some $\xi$ such that
\begin{equation*}
w(\xi)=\min_{\Sigma_{\la}}w <0.
\end{equation*}
At the negative minimum point of $w$, we have
\begin{eqnarray}\nonumber
  \lap U(\xi)&=& \lap w(\xi)\, \phi(\xi)+2 \nabla w (\xi)\cdot \nabla \phi(\xi)
  +w(\xi)\, \lap \phi(\xi) \\\nonumber
   &\geq& - w(\xi)\, \frac{\phi(\xi)}{l^2}\\\nonumber
     &\geq&  - U(\bar{x})\,\frac{\phi(\xi)}{\phi(\bar{x})\,l^2}.
\end{eqnarray}
On the other hand, by (\ref{2016541}) and (\ref{2016527}), we know
\begin{eqnarray}\nonumber
\lap U(\xi) &\leq & c_1(\xi)\,V(\xi) \\\label{20165123}
   &\leq&   -c_1(\xi)\,c_2(x^0)\,U(x^0)\, l^\alpha.
\end{eqnarray}
Combining
(\ref{20165122}) and (\ref{20165123}), it yields
$$l^{2+\alpha}\:\frac{\phi(\bar{x})\,c_1(\xi)\,c_2(x^0)}{\phi(\xi)}\geq 1.$$
The inequality above will certainly fail to hold for
$l$ sufficiently small.
This shows that $V(x)$ must be non-negative. Applying the \emph{maximum principle}
to (\ref{2016541}), it yields
$$U(x)\geq0, \quad x \in \Omega.$$

Next we  argue (\ref{20165124}) by a contradiction.

Suppose for some $\eta \in \Omega$, $U(\eta)=0$. Then $\eta$ is the minimum point of $U$. Thus
$$0\geq-\lap U(\eta)\geq c_1(\eta) V(\eta).$$
Meanwhile,
$$c_1(\eta) V(\eta)\leq0.$$
Hence
$$V(\eta)=0= \min_{\Sigma_{\la}} V,$$
and
\begin{eqnarray}
&&(-\lap)^{\alpha/2}V(\eta)\nonumber\\
&=&C_{n,\alpha}PV\int_{\mathbb{R}^n}
\frac{-V(y)}
{|\eta-y|^{n+\alpha}}dy\nonumber\\
&=&C_{n,\alpha}PV\int_{\Sigma_{\la}} \frac{-V(y)}
{|\eta-y|^{n+\alpha}}dy+\int_{\Sigma_{\la}} \frac{-V(y^{\la})}
{|\eta-y^{\la}|^{n+\alpha}}dy\nonumber\\\label{20165125}
&=&C_{n,\alpha}PV\int_{\Sigma_{\la}} \bigg(\frac{1}
{|x^0-y^{\la}|^{n+\alpha}}-\frac{1}
{|x^0-y|^{n+\alpha}} \bigg) \,V(y)\,dy.
\end{eqnarray}
If $V(x)\not\equiv0$, then (\ref{20165125}) implies that
$$(-\lap)^{\alpha/2}V(\eta)<0.$$
Together with (\ref{2016541}), it shows that
$$U(\eta)<0.$$
This is a contradiction with  (\ref{20165126}). Hence $ V(x)$ must be identically 0 in $\Sigma_{\la}$. Since $$V(x^{\la})=-V(x),  x \in \Sigma_{\la},$$
it shows that
$$V(x)\equiv0, \quad x \in R^n.$$
 Again with (\ref{2016541}), one can easily deduce that
$U(x)\leq 0, x \in \Sigma_{\la}$. Since we already know that
$U(x)\geq 0, x \in \Sigma_{\la}$,
it must hold that
$U(x)= 0, x \in \Sigma_{\la}$.
Together with $U(x^{\la})=U(x)$, we arrive at
$$U(x)\equiv 0, x \in R^n.$$

If there exists a $\xi \in \Omega$ such that $V(\xi)=0$. Then
from (\ref{2016527}) and  (\ref{2016541}) it follows
$$-c_1(\xi)V(\xi)\leq(-\lap)^{\alpha/2} V(\xi)\leq \frac{V(\xi)}{l^\alpha}<0.$$
It thus implies that some $\eta \in \Omega$, $U(\eta)=0$. The rest is the same as
the previous argument. It thus proves (\ref{20165124}).

\subsection{The Iteration Method}

\begin{thm}
Let $s=1+\alpha/2$ for $0<\alpha<2$. Assume that $u \in C_{loc}^{3,1}(B_1(0))$
is a positive solution of
\begin{equation}\label{2016525}
\left\{\begin{array}{ll}
(-\lap)^s u=f(u),  &x \in B_1,\\
u=\lap u=0, &x \in R^{n}\backslash B_1,
\end{array}
\right.
\end{equation}
where $f(t)$ is Lipschitz continuous and increasing in $t$. Then $u$ must be symmetric about the origin.
\end{thm}

\textbf{Proof. }
Let $-\lap u=v$. Then (\ref{2016525}) can be split into two equations:
\begin{equation}
\left\{\begin{array}{ll}
(-\lap)^{\alpha/2}v=f(u),  &x \in B_1,\\
v=0, &x \in R^{n}\backslash B_1,
\end{array}
\right.
\end{equation}
and
\begin{equation}
\left\{\begin{array}{ll}
-\lap u=v,  &x \in B_1,\\
u=0, &x \in R^{n}\backslash B_1.
\end{array}
\right.
\end{equation}

We carry out the proof via the \emph{method of moving planes}.

We first choose any direction to be the $x_1$-direction and let
$$T_\lambda=\{x\in\mathbb{R}^n \mid x_1=\lambda\}, \; \;
\Sigma_\lambda=\{x\in \mathbb{R}^n  \mid x_1<\lambda\},$$
and $$x^\lambda=\{(2\lambda-x_1, x^\prime) \mid x=(x_1, x^\prime)\in \mathbb{R}^n \}$$
be the reflection of $x$ about the plane $T_\lambda$.

Let $$u_\lambda(x)=u(x^\lambda)$$ and
$$U_\lambda(x)=u_\lambda(x)-u(x), \qquad
V_\lambda(x)=v_\lambda(x)-v(x).$$
Then
\begin{equation}\label{2016544}
\left\{\begin{array}{ll}
-\lap U_{\la}(x)=V_{\la}(x), &x \in \Sigma_\lambda,\\
(-\lap)^{\alpha/2} V_{\la}(x)=f(u_{\la})-f(u),
&x \in \Sigma_\lambda\cap B_1(0),\\
U_{\la}(x)\geq0, &x \in  \partial (\Sigma_\lambda \cap B_1(0)),\\
V_{\la}(x)\geq0,& x \in \Sigma_{\la}\backslash (\Sigma_\lambda \cap B_1(0)).
\end{array}
\right.
\end{equation}

\emph{Step 1. Moving the plane $T_{\la}$ from $-1$ to the right.}

For $\la$ near $-1$, we claim that
\be
U_\lambda(x)\geq 0, \quad x \in \Sigma_\lambda.
\ee

Notice that
$$f(u_{\la})-f(u)=\frac{f(u_{\la})-f(u)}{u_{\la}(x)-u(x)} U_{\la}(x),$$
this allows us to obtain desired result
by applying Lemma \ref{2016541} (\textit{narrow region principle})
 to (\ref{2016544}) with $c_1(x)=-1$ and
$$c_2(x)=-\frac{f(u_{\la})-f(u)}{u_{\la}(x)-u(x)}.$$

\emph{Step 2. Keep moving the plane $T_{\la}$ until the limiting position}
$$\la_0=\sup\{\la \leq 0\mid U_\mu(x), V_\mu(x)\geq 0, x \in \Sigma_\mu, \forall \mu < \la \}.$$

We claim that
\be\label{2016545}
\la_0=0.
\ee
If not, then for $\la_0<0$ we can show that
\be\label{2016546}
U_{\la_0}(x)\equiv0, \quad x \in \Sigma_{\lambda_0}.
\ee
We postpone the proof of (\ref{2016546}) for the moment. From
(\ref{2016545}), we know that
$$
0=u_{\la_0}(x)=u(x)>0, \quad x=1+\la_0.
$$
A contradiction. This proves (\ref{2016545}).

To prove (\ref{2016546}), we argue by contradiction. If
(\ref{2016546}) does not hold, then by (\ref{20165124}) (the strong maximum principle), we have
$$U_{\la_0}(x)>0, \quad x \in \Sigma_{\la_0} \cap B_1(0).$$
This enables us to keep moving the plane $T_{\la_0}$ to the right. Precisely
speaking, for $\varepsilon>0$ small such that $\la_0+\varepsilon<0$, it holds
\begin{equation}\label{2016547}
 U_{\la}(x)\geq0, \quad x \in \Sigma_{\la} \cap B_1(0) , \quad \la \in [\la_0,\la_0+\varepsilon).
\end{equation}
The inequality above contradicts the definition of $\la_0$. Therefore (\ref{2016546})
must be true.

Now we show (\ref{2016547}) by constructing a narrow region. For $\delta>0$ small,
there exists a constant $C$ such that
$$
   U_{\la_0}(x)\geq C>0, \quad x \in \Sigma_{\la-\delta} \cap B_1(0).
$$
By the continuity of $U$ in $\la$, we have
\begin{equation}\label{2016548}
  U_{\la}(x)\geq 0, \quad x \in \Sigma_{\la-\delta} \cap B_1(0).
\end{equation}
It's easy to see that
$$
  U_{\la}(x)\geq 0, \quad x \in
  \partial((\Sigma_{\la}\backslash\Sigma_{\la-\delta})\cap B_1).
$$
Since $\Sigma_{\la}\backslash\Sigma_{\la-\delta}$ is a narrow region, with Lemma \ref{2016541} we derive that
$$
  U_{\la}(x)\geq 0, \quad x \in
 (\Sigma_{\la}\backslash\Sigma_{\la-\delta})\cap B_1.
$$
Together with (\ref{2016548}), we arrive at (\ref{2016547}). This completes with
proof of (\ref{2016545}).

Now we have
$$U_0(x)\geq 0, \quad x \in \Sigma_0.$$
Similarly, one can move $T_{\la}$ from near 1 to the left and show that
$$U_0(x)\leq 0, \quad x \in \Sigma_0.$$
Thus
$$U_0(x)\equiv 0, \quad x \in \Sigma_0.$$
Since the direction of the $x_1$-axis is arbitrary, we have actually proved that
$u$ is symmetric about the origin.

Authors' Addresses and Emails；

Ran Zhuo

Department of Mathematical Sciences,

Huanghuai University, Zhumadian, China, 463000

zhuoran1986@126.com

\smallskip

Yan Li

Department of Mathematical Sciences,

Yeshiva University, New York, NY, 10033

yali3@mail.yu.edu

\end{document}